\documentclass[12pt, leqno]{amsart}

\usepackage{shuffle}

\usepackage{lscape}
\usepackage{graphicx}
\usepackage{amssymb,amsmath,amsthm,amscd}
\usepackage{mathrsfs}
\usepackage{enumerate}
\usepackage[usenames,dvipsnames]{color}
\usepackage[colorlinks=true, pdfstartview=FitV,
 linkcolor=blue,citecolor=blue,urlcolor=blue]{hyperref}

\usepackage[all]{xy}

\allowdisplaybreaks[4]
\usepackage{oldgerm}
\newcommand{\nc}{\newcommand}
\numberwithin{equation}{section}

\newenvironment{red}{\relax\color{red}}{\relax}
\newenvironment{jaune}{\relax\color{Orchid}}{\hspace*{.5ex}\relax}

\newcommand{\bj}{\begin{jaune}}
\newcommand{\ej}{\end{jaune}}

\newcommand{\ber}{\begin{red}}
\newcommand{\er}{\end{red}}

\theoremstyle{plain}
\newtheorem{lemma}{Lemma}[section]
\newtheorem{prop}[lemma]{Proposition}
\newtheorem{theorem}[lemma]{Theorem}

\newcommand{\Prop}{\begin{prop}}
\newcommand{\enprop}{\end{prop}}
\newcommand{\Lemma}{\begin{lemma}}
\newcommand{\enlemma}{\end{lemma}}
\newcommand{\Th}{\begin{theorem}}
\newcommand{\enth}{\end{theorem}}
\newtheorem{corollary}[lemma]{Corollary}
\newcommand{\Cor}{\begin{corollary}}
\newcommand{\encor}{\end{corollary}}
\newtheorem{definition}[lemma]{Definition}
\newtheorem*{conjecture}{Conjecture}
\newcommand{\Def}{\begin{definition}}
\newcommand{\edf}{\end{definition}}
\newtheorem{sublemma}[lemma]{Sublemma}
\newcommand{\Sublemma}{\begin{sublemma}}
\newcommand{\ensub}{\end{sublemma}}

\theoremstyle{definition}
\newtheorem{remark}[lemma]{Remark}

\newtheorem{Convention}[lemma]{Convention}
\newcommand{\Conv}{\begin{Convention}}
\newcommand{\enconv}{\end{Convention}}
\nc{\Conj}{\begin{conjecture}}
\nc{\enconj}{\end{conjecture}}
\nc{\Rem}{\begin{remark}}
\nc{\enrem}{\end{remark}}
\newcommand{\C}{{\mathbb C}}

\newcommand{\Z}{{\mathbb Z}}
\newcommand{\B}{{\mathbf{B}}}
\newcommand{\A}{{\mathbf A}}

\newcommand{\D}{\mathscr{D}}

\newcommand{\one}{{\bf{1}}}
\newcommand{\seteq}{\mathbin{:=}}

\newcommand{\hd}{{\operatorname{hd}}}

\newcommand{\g}{{\mathfrak{g}}}

\newcommand{\Hom}{\operatorname{Hom}}
\newcommand{\End}{\operatorname{End}}

\newcommand{\isoto}[1][]{\mathop{\xrightarrow%
[{\raisebox{.3ex}[0ex][.3ex]{$\scriptstyle{#1}$}}]%
{{\raisebox{-.6ex}[0ex][-.6ex]{$\mspace{2mu}\sim\mspace{2mu}$}}}}}

\newcommand{\M}{{\mathscr M}}

\newcommand{\eq}{\begin{eqnarray}}
\newcommand{\eneq}{\end{eqnarray}}

\newcommand{\hs}{\hspace*}

\newcommand{\To}[1][{\hs{2ex}}]{\xrightarrow{\,#1\,}}

\newcommand{\eqn}{\begin{eqnarray*}}
\newcommand{\eneqn}{\end{eqnarray*}}
\newcommand{\on}{\operatorname}

\newcommand{\bna}{\be[{\rm(a)}]}

\newcommand{\soplus}{\mathop{\mbox{\normalsize$\bigoplus$}}\limits}

\newcommand{\ba}{\begin{array}}
\newcommand{\ea}{\end{array}}

\newcommand{\indlim}{\varinjlim\limits}
\newcommand{\set}[2]{\left\{#1 \mid #2 \right\}}

\newcommand{\eqsub}{\begin{subequations}\begin{eqnarray}}
\newcommand{\eneqsub}{\end{eqnarray}\end{subequations}}

\newcommand{\ol}{\overline}

\nc{\la}{\lambda}
\nc{\lam}{\lambda}
\nc{\U}[1][\g]{U_q(#1)}
\nc{\te}{\tilde{e}}
\nc{\tei}{\tilde{e}_i}
\nc{\tf}{\tilde{f}}
\nc{\tfi}{\tilde{f}_i}
\nc{\tU}{\widetilde U_q(\g)}
\nc{\tE}{\tilde{E}}
\nc{\tF}{\widetilde{\F}}
\nc{\tK}{\widetilde{K}}

\nc{\tk}{\tilde{k}}
\nc{\tkone}{\tk_{\ol{1}}}
\nc{\teone}{\tilde{e}_{\ol{1}}}
\nc{\tfone}{\tilde{f}_{\ol{1}}}

\nc{\teibar}{\tilde{e}_{\ol{i}}} \nc{\tfibar}{\tilde{f}_{\ol{i}}}
\nc{\tki}{{\tk}_{\ol {i}}}

\nc{\BZ}{{\mathbb{Z}}}
\nc{\al}{\alpha}
\nc{\qs}{{q}}
\nc{\lan}{\langle}
\nc{\ran}{\rangle}
\nc{\re}{{\mathrm{re}}}
\nc{\wt}{\operatorname{wt}}
\nc{\ch}{\operatorname{ch}}
\nc{\Um}[1][\g]{U^-_q(#1)}
\nc{\Ue}{U^+_q(\g)}
\nc{\eps}{\varepsilon}
\nc{\vphi}{\varphi}
\nc{\sphi}{\varphi^*}
\nc{\seps}{\varepsilon^*}

\nc{\nn}{\nonumber}

\nc{\vp}{\varpi}
\nc{\cls}{{\operatorname{cl}}}
\nc{\Wt}{{\operatorname{Wt}}}
\nc{\Us}{U'_q(\g)}
\nc{\La}{\Lambda}
\nc{\tLa}{\widetilde\Lambda}
\nc{\ro}{{\rm(}}
\nc{\rf}{{\rm)}}
\nc{\norm}{{\mathrm{norm}}}
\nc{\qbox}{\quad\mbox}
\nc{\braid}{{\mathfrak{B}}}
\nc{\Ad}{\operatorname{Ad}}
\nc{\Aut}{\operatorname{Aut}}
\nc{\dt}[1]{\tilde{\tilde #1}}
\nc{\Sn}{S^{{\mathrm{norm}}}}
\nc{\aff}{{\rm{aff}}}
\nc{\rk}{{\mathrm{rk}}}
\nc{\tP}{\widetilde{P}}
\nc{\tW}{\widetilde{W}}
\nc{\Dyn}{\mathrm{Dyn}}
\nc{\tD}{\widetilde{\Delta}}
\nc{\height}[1]{{\operatorname{ht}}(#1)}
\nc{\bl}{\bigl(}
\nc{\br}{\bigr)}
\nc{\Hecke}{\mathrm{H}}
\nc{\HA}{\Hecke^{\mathrm{A}}}
\nc{\HB}{\Hecke^{\mathrm{B}}}
\newcommand{\scbul}{{\,\raise1pt\hbox{$\scriptscriptstyle\bullet$}\,}}
\nc{\vac}{{\phi}}
\nc{\Bt}{\B_\theta(\g)}
\nc{\be}{\begin{enumerate}}
\nc{\ee}{\end{enumerate}}
\nc{\low}{{\mathrm{low}}}
\nc{\upper}{{\mathrm{up}}}
\nc{\Zodd}{\Z_{\mathrm{odd}}}
\nc{\Ft}[1][n]{\mathbb{P}\mathrm{ol}_{#1}}
\nc{\Ftf}[1][n]{\widetilde{\mathbb{P}\mathrm{ol}}_{#1}}
\nc{\KA}{\on{K}^{\mathrm{A}}}
\nc{\KB}{\on{K}^{\mathrm{B}}}
\nc{\Res}{\on{Res}}
\nc{\Fc}[1][{n,m}]{\mathbf{F}_{#1}}
\nc{\tphi}{\tilde{\varphi}}
\nc{\CO}{\mathscr{O}}
\nc{\inte}{\mathrm{int}}
\nc{\Oint}{\mathcal{O}^{\ge0}_{\inte}}
\nc{\vs}{\vspace*}
\nc{\tLt}{\widetilde{L}}
\nc{\tL}{\widetilde{\Lambda}}
\nc{\tu}{\tilde{u}}
\nc{\noi}{\noindent}
\nc{\heigh}{\mathfrak{t}}
\nc{\lowest}{\mathfrak{l}}
\nc{\rootl}{\mathsf{Q}}
\nc{\cl}{\colon}
\nc{\uqpg}{U'_q(\mathfrak g)}
\nc{\uq}{\uqpg}
\nc{\Oh}{\widehat{\mathcal{O}}}
\nc{\pn}{p_{\mathfrak{n}}}

\nc{\KLR}{KLR algebra}
\nc{\KLRs}{KLR algebras}
\nc{\cor}{\mathbf{k}}
\nc{\cora}{{\cor(A)}}
\nc{\haut}{\mathrm{ht}}
\nc{\tens}{\mathop\otimes}
\nc{\gmod}{\mbox{-$\mathrm{gmod}$}}
\nc{\gMod}{\mbox{-$\mathrm{gMod}$}}
\nc{\proj}{\mbox{-$\mathrm{proj}$}}
\nc{\gproj}{\mbox{-$\mathrm{gproj}$}}
\nc{\smod}{\mbox{-$\mathrm{mod}$}}
\nc{\Mod}{\mbox{-$\mathrm{Mod}$}}
\nc{\h}{\mathfrak h}
\nc{\Rnorm}{R^{\rm{norm}}}
\nc{\Runiv}{R^{\rm{univ}}}
\nc{\Rren}{R^{\rm{ren}}}

\nc{\Vhat}{\widehat{V}}
\nc{\F}{\mathcal{F}}

\def\T{{\mathcal T}}

\nc{\fd}[1][A]{\on{\mathrm{flat.dim}_{#1}}}
\nc{\bP}{{\mathbb{P}}}
\nc{\bPh}{\widehat{\mathbb{P}}}
\nc{\bK}[1][{n}]{\widehat{\mathbb{K}}_{#1}}
\nc{\bV}[1][{n}]{\widehat{V}^{\otimes{#1}}}
\nc{\bVK}[1][{n}]{\widehat{V}^{\otimes{#1}}_{\widehat{\mathbb{K}}}}
\nc{\hV}{\widehat{V}}
\nc{\opp}{\mathrm{opp}}
\nc{\col}{\colon}
\nc{\bnum}{\be[{\rm(i)}]}
\nc{\oep}{\epsilon}
\nc{\qtext}{\quad\text}
\nc{\qtextq}[1]{\quad\text{#1}\quad}
\nc{\longtwoheadrightarrow}[1][]{\xymatrix{\ar@{->>}[r]^-{{#1}}&}}
\nc{\epiTo}[1][]{\longtwoheadrightarrow[{#1}]}
\nc{\epito}{\twoheadrightarrow}
\nc{\monoTo}[1][]{\xymatrix{\ar@{>->}[r]^-{{#1}}&}}
\nc{\sym}{\mathfrak{S}}
\nc{\inp}[1]{{({#1})_{\mathrm{n}}}}
\nc{\rtl}{\rootl}
\nc{\wtd}{\widetilde}
\nc{\etens}{\boxtimes}
\nc{\ds}[1]{\mathrm{d}(#1)}
\nc{\rmat}[1]{{\mathbf{r}}_%
{\mspace{-2mu}\raisebox{-.6ex}{${\scriptstyle{#1}}$}}}
\nc{\rmats}[1]{{\mathbf{r}}_%
{\mspace{-2mu}\raisebox{-.6ex}{${\scriptscriptstyle{#1}}$}}}
\nc{\shc}{\mathcal{C}}
\nc{\shs}{\mathcal{S}}
\nc{\Fct}{{\on{Fct}}}
\nc{\tC}{\widetilde{\shc}}
\nc{\Zp}{\Z_{\ge0}}
\nc{\tPhi}{\widetilde{\Phi}}
\nc{\tT}{{\widetilde{\T}}}
\nc{\Ob}{\on{Ob}}
\nc{\bwr}{\mbox{\large$\wr$}}
\nc{\Img}{\on{Im}}
\nc{\Ab}{\mathcal{A}^{\mathrm{big}}}
\nc{\Sb}{\mathcal{S}^{\mathrm{big}}}
\nc{\As}{\mathcal{A}}
\nc{\Ss}{\mathcal{S}}
\nc{\ntens}{\widetilde{\otimes}}
\nc{\hR}{\widehat{R}}
\nc{\nconv}{\mathop{\mbox{\large $\odot$}}}
\nc{\snconv}{\mbox{\scriptsize$\odot$}}
\nc{\ts}{\tilde{s}}
\nc{\sho}{\mathcal{O}}
\nc{\bc}{\begin{cases}}
\nc{\ec}{\end{cases}}
\nc{\slnh}{{\widehat{\mathfrak{sl}}_N}}
\nc{\UA}{U_q'(\slnh)}
\nc{\KR}{R_K}
\nc{\cQ}{\mathcal{Q}}
\nc{\Irr}{\mathcal{I}rr}
\nc{\tQ}{\widetilde{\cQ}}
\nc{\bs}{\mathbf{s}}
\nc{\bL}{\mathbb{L}}
\nc{\tg}{\tilde{g}}

\nc{\conv}{\mathbin{\mbox{\large $\circ$}}}
\nc{\shconv}{\mathbin{\large\diamond}}
\nc{\sconv}{\mathbin{\large\Delta}}
\nc{\hconv}{\mathbin{\nabla}}

\nc{\Rm}{R^{\mathrm{ren}}}

\nc{\bQ}{\ol{Q}}

\nc{\de}{\on{\textfrak{d}}}

\nc{\xmono}{\ar@{>->}}
\nc{\xepi}{\ar@{->>}}
\nc{\db}[1]{\raisebox{-.5ex}[2ex][1.8ex]{$#1$}}
\nc{\wb}[1]{\mbox{$\rule[-1.1ex]{0ex}{2ex}#1$}}
\nc{\univ}{\mathrm{univ}}
\nc{\rM}{{}^*\mspace{-2mu}M}
\nc{\lM}{M^*}
\nc{\uqm}{\uq\smod}
\nc{\tR}{\widetilde{R}_{\gamma,\beta}}
\nc{\tx}{\tilde{x}}
\nc{\bi}{\mathbf{i}}
\nc{\ttau}{\widetilde{\tau}}

\nc{\tEnd}{\on{\widetilde{E}nd}}
\nc{\tHom}{\on{\widetilde{H}om}}

\nc{\K}{{J}}
\nc{\Kex}{{\K}_{\mathrm{ex}}}
\nc{\Kfr}{{\K}_{\mathrm{f\mspace{.01mu}r}}}
\nc{\coro}{\cor}
\nc{\tB}{\widetilde{B}}
\nc{\seed}{\mathscr{S}}

\nc{\up}{\mathrm{up}}
\nc{\bfa}{\mathbf{a}}

\newlength{\mylength}
\nc{\ov}[1]{\overline{#1}}
\nc{\Wlmj}[3]{\W_{#2,#3}^{(#1)}}
\nc{\Mkl}[2]{\M_\ttww(#1,#2)}
\nc{\mqs}{(-q^2)}
\nc{\Cquiver}{\upsigma}

\nc{\mut}[1]{{\mu}_{\mspace{-2mu}\raisebox{-.5ex}{${\scriptstyle{#1}}$}}}

\nc{\Kt}{\mathcal K_t}
\nc{\KT}{\mathbb{K}_t}
\nc{\yim}{y_{i,m}}
\nc{\yjm}{y_{j,m}}
\nc{\yjp}{y_{j,p}}
\nc{\yimp}{y_{i,m+1}}
\nc{\yjmp}{y_{j,m+1}}

\nc{\Refl}{\mathscr{S}}
\nc{\Reflinv}{{\Refl}^{-1}}

\nc{\catC}{\mathscr C}
\nc{\catA}{\mathcal A}
\nc{\shift}{{\mathrm T}}
\nc{\rE}{ \mathsf{E} }
\nc{\rW}{ \mathcal{W} }
\nc{\rES}{ \mathcal{E} }

\nc{\brd}{\sigma} 
\nc{\into}{\xymatrix@C=3ex{{}\ar@{^{(}->}[r]&{}}}
\nc{\dual}{\D}
\nc{\cat}[1][{\g}]{\catC_{#1}^0}
\nc{\qt}[1]{[{#1}]_t}

\title[{\scalebox{.9}{Braid group action on the module category of 
quantum affine algebras}}]
{Braid group action on the module category of quantum affine algebras}

\author[M. Kashiwara]{Masaki Kashiwara}
\thanks{The research of M.\ Kashiwara
was supported by Grant-in-Aid for Scientific Research (B)
15H03608, Japan Society for the Promotion of Science.}
\address[M. Kashiwara]{
Kyoto University Institute for Advanced Study,
Research Institute for Mathematical Sciences, Kyoto University,
Kyoto 606-8502, Japan \& Korea Institute for Advanced Study, Seoul 02455, Korea }
\email[M. Kashiwara]{masaki@kurims.kyoto-u.ac.jp}

\author[M. Kim]{Myungho Kim}
\address[M. Kim]{Department of Mathematics, Kyung Hee University, Seoul 02447, Korea}
\email[M. Kim]{mkim@khu.ac.kr}
\thanks{The research of M.\ Kim was supported by the National Research Foundation of
Korea(NRF) Grant funded by the Korea government(MSIP) (NRF-2017R1C1B2007824).}

\author[S.-j. Oh]{Se-jin Oh}
\thanks{ The research of S.-j.\ Oh was supported by the Ministry of Education of the Republic of Korea and the National Research Foundation of Korea (NRF-2019R1A2C4069647).}
\address[S.-j. Oh]{Department of Mathematics, Ewha Womans University, Seoul 03760, Korea}
\email[S.-j. Oh]{sejin092@gmail.com}

\author[E. Park]{Euiyong Park}
\address[E. Park]{Department of Mathematics, University of Seoul, Seoul 02504, Korea}
\email[E. Park]{epark@uos.ac.kr}


\date{April 10, 2020}

\begin{document}
\maketitle
\begin{abstract}
Let $\g_0$ be a simple Lie algebra of type ADE and let $\uqpg$ be the corresponding untwisted quantum affine algebra. We show that  there exists an action of the braid group $B(\g_0)$ on the quantum Grothendieck ring $\Kt(\g)$ of Hernandez-Leclerc's category $\catC_\g^0$. 
Focused on the case of type $A_{N-1}$,
we construct a family of monoidal autofunctors  $\{\Refl_i\}_{i\in \Z}$ on a localization  $\T_N$ of the category of finite-dimensional graded modules over the quiver Hecke algebra of type $A_{\infty}$.  Under an isomorphism between the Grothendieck ring $K(\T_N)$ of $\T_N$  and the quantum Grothendieck ring  $\Kt({A^{(1)}_{N-1}})$, the functors $\{\Refl_i\}_{1\le i\le N-1}$ recover the  action of the braid group $B(A_{N-1})$. We investigate further properties of these functors.
\end{abstract}

\section{Introduction}
The monoidal category $\catC_\g$ of  finite-dimensional representations of a quantum affine algebra $\uqpg$ has been extensively investigated because of its rich structure. 
Among various approaches,  Nakajima (\cite{Nak04}), Varagnolo-Vasserot (\cite{VV03}), and Hernandez (\cite{Her04})  studied $t$-deformations of the Grothendieck ring of $\catC_\g$. 
These $t$-deformations are interesting, because they provide a way to calculate the $q$-character of simple representations. 
There is  a full subcategory $\catC^0_\g$ of $\catC_\g$, introduced by Hernandez and Leclerc in \cite{HL10},   which contains an essential information 
on $\catC_\g$ but  has a smaller set of the classes of simple modules. 
The Grothendieck ring of  $\catC^0_\g$ is isomorphic to the polynomial ring in countably many variables, while that of $\catC_\g$ is the one in  uncountably many variables.
For the cases where $\g$ is one  of  untwisted ADE types, a $t$-deformation  $\Kt(\g)$ of the Grothendieck ring of $\catC^0_\g$, called the \emph{quantum Grothendieck ring},  
was investigated from a ring theoretic point of view in \cite{HL13}.
It turns out that the $\C(t^{1/2})$-algebra $\Kt(\g)$ has an interesting presentation: there is  a set of generators consisting of a countable infinite number of copies of Drinfeld-Jimbo generators of a half  of the quantum group $U_t (\g_0)$,  and they satisfy the quantum Serre relations in a copy, $t$-boson relations between adjacent copies,  and $t$-commutation relations between  non-adjacent copies. 
This presentation reflects the following feature of the category   $\catC_\g^{0}$: 
for each choice of a Dynkin quiver $Q$ with an additional data,  they defined
a monoidal subcategory $\catC_Q$ of $\catC_\g^{0}$ such that  the quantum Grothendieck ring of $\catC_Q$ is isomorphic to the half  $U^-_t (\g_0)$ of the quantum group $U_t (\g_0)$, and  all the fundamental representations in  $\catC_\g^{0}$ can be obtained from those in $\catC_Q$ by taking 
functors $\D^m$ ($m\in\Z$).
Here $\D$ is the contravariant
functor taking the right dual.

One of  main results of this paper is that there exists an action of the braid group $B(\g_0)$ of type $\g_0$ on the quantum Grothendieck ring $\Kt(\g)$ (Theorem~\ref{thm: braid group action}).
Since we give the action explicitly, the braid relations  can be checked by the presentation of $\Kt(\g)$. 
Recall that the blocks of the category $\catC_\g^0$ is parameterized by the root lattice of $\g_0$ and the tensor product is compatible with the addition on the root lattice (\cite{KKOP20}).
It turns out that the action of  the generators $\brd_i$ of $B(\g_0)$ on $\Kt(\g)$ correspond to the reflections with respect to the simple roots $\al_i$ on the root lattice. 
Indeed, the action of $\brd_i$ is related with Saito's reflection functor
as seen in  Theorem~\ref{Saito}.

We conjecture that the braid group action can be lifted to the action on the monoidal category $\catC_\g^0$.
We show that it is the case when $\g$ is of  type $A_{N-1}^{(1)}$. 
A key point of view is the use of a rigid monoidal category $\T_N$ which is constructed out of the category $\catA$ of finite-dimensional graded modules over the quiver Hecke algebra $R^{A_\infty}$ of type $A_\infty$ (\cite{KKK18A}).
It is a certain localization of $\catA$ and  
there is a monoidal functor $\F_N$ from $\T_N$ to $\catC^0_{A_{N-1}^{(1)}}$ which sends simple objects to simple objects. Moreover this functor induces an isomorphism between the Grothendieck ring  
$K(\T_N)$ and the quantum Grothendieck ring $\Kt(A_{N-1}^{(1)})$.
It is summarized by the diagram
$$K(\T_N)\isoto \Kt(A_{N-1}^{(1)})\To[\ {t=1}\ ] K(\catC^0_{A_{N-1}^{(1)}}).$$
Hence  the category $\T_N$ can be understood as a \emph{graded lift} of $\catC^0_{A_{N-1}^{(1)}}$ as a rigid monoidal category.

We  show that there is a family of monoidal autofunctors  $\{\Refl_i\}_{1\le i \le N-1}$ on the category  $\T_N$ which recover the action of the braid group $B(A_{N-1})$ under the  isomorphism between $K(\T_N)$  and $\Kt({A^{(1)}_{N-1}})$ (Theorem \ref{thm: reflection functors}, Theorem \ref{thm:functors recover the action}).
There is a  general procedure, developed in \cite{KP},  to construct monoidal functors between the  categories of modules over quiver Hecke algebras, and a similar procedure can be applied for the category $\T_N$. 
This is a main advantage in working on the category $\T_N$ rather than the category $\catC^{0}_{A_{N-1}}$.

Finally we  provide several consequences of 
the existence of such functors $\Refl_i$. 
For a simple object $L$ which belongs to an orbit of $L(i)$ 
for some $i$ under the action $B(A_{N-1})$, 
one can define an automorphism $s_L$ which has similar properties
with the automorphisms $s_i$ (Theorem \ref{thm: reflection by root modules}). 
Moreover $s_{L(i)}$ coincides with $s_i$.

\smallskip
This paper is an announcement whose details will appear elsewhere.

\section{Braid group action on the quantum Grothendieck rings}

\newcommand{\modt}{ \; \mathrm{mod}\; 2  }
Let $\g_0$ be a finite-dimensional simple Lie algebra of simply-laced type with a Cartan matrix $\A=(a_{ij})_{i,j\in I_0}$, 
$\g$ the untwisted affine Kac-Moody algebra associated with $\g_0$,
 and $\uqpg$ the quantum affine algebra associated with $\g$.
 We take the algebraic closure of $\C(q)$ inside $\bigcup_{m >0}\C((q^{1/m}))$
 as the base field  $\cor$ for $\uqpg$. 
 Let $\catC_\g$ be the category of finite-dimensional integrable modules over $\uqpg$.
 There is a family $\set{V(\varpi_i)_c}{i\in I_0,c\in \cor^\times}$ in $\catC_\g$, called the \emph{fundamental representations}.

Following \cite{HL10}, we denote by 
$\catC_\g^0$ the smallest full subcategory of the category $\catC_\g$ 
which is stable under taking subquotients, extensions,  
tensor products and contains
$$\{V(\varpi_i)_{(-q)^{p}} \ | i \in I_0,  \ p \equiv d(1,i) \modt \},$$
where $d(i,j)$ is the distance between the vertices $i$ and $j$ in the Dynkin diagram of $\g_0$. Here $1\in I_0$ is an arbitrarily chosen element.
Then the  complexified Grothendieck ring $\C\tens_\Z K(\catC_\g^0)$ of 
$\catC_\g^0$ has a $t$-deformation $\Kt(\g)$, called the \emph{quantum Grothendieck ring} of $\catC_\g^0$.  
To each simple module $S$ in $\cat$, we can associate
an element $\qt{S}$ of $\Kt(\g)$ and we have
$\Kt(\g)=\soplus_S\C(t^{1/2})\qt{S}$. Here $S$ ranges over the set of the
isomorphism classes of simple modules in $\cat$.

Let $Q$ be a Dynkin quiver with type $\g_0$, 
and let $\phi_Q$ be a height function, i.e.,
it associates an integer $\phi_Q(i)$ to each vertex $i$ of $Q$ such that
$\phi_Q(i)=\phi_Q(j)+1$ if $i\to j$.
We assume further that $\phi_Q(1)\in2\Z$.
A pair $Q=(Q,\phi_Q)$ is called a Q-data. 

For a sink $i$ of $Q$, let $s_iQ\seteq(s_iQ,\phi_{s_iQ})$ be the Q-data consisting 
of the Dynkin quiver $s_iQ$ obtained from $Q$ by reversing the arrows
of $Q$ adjacent to $i$ and the height function $\phi_{s_iQ}$ of $s_iQ$
given by
$\phi_{s_iQ}(j)=\phi_Q(j)+2\delta_{i,j}$.

To a Q-data $Q$,
Hernandez-Leclerc (\cite{HL13}) associated a full
monoidal subcategory $\catC_Q$ of $\catC^0_\g$, and 
a monoidal functor
$\F_Q\col R_{\g_0}\smod\to \catC_Q$ is constructed
in \cite{KKK15B,Fujita18}, and Fujita (\cite{Fujita, Fujita18}) 
proved that $\F_Q$ is an equivalence of categories.
Here, $R_{\g_0}\smod$ is the monoidal category of finite-dimensional modules
(with nilpotent actions of the generators $x_k$) over the quiver Hecke algebra $R_{\g_0}$ associated with $\g_0$.
Note that $\F_Q(L(i))$ is a fundamental module for any $i\in I_0$, where
$L(i)\in R_{\g_0}\smod$ is the simple module associated with $i\in I_0$.

Then, for a Q-data $Q$, we have
an embedding of $\Z[t^{\pm1}]$-algebras
$$j_Q\col K(R_{\g_0}\gmod)\into\Kt(\g)$$
induced by $\F_Q$. %

Let $\KT(\g_0)$ be 
the $\C(t^{1/2})$-algebra 
 generated by 
$\set{y_{i,m}}{i \in I_0, m \in \Z}$
with the defining relations:
\bna
\item for $m\in \Z$,
\eqn 
\begin{aligned}
&\hs{-1em}\yim \yjm=\yjm\yim   \quad \text{if } \ a_{ij}=0, \\
&\hs{-1em}\yim^2 \yjm -(t+t^{-1})\yim \yjm \yim  + \yjm \yim^2  =0  
\quad\text{if $a_{ij}=-1$, }
\end{aligned}
\eneqn
\item for $m \in \Z$ and $i,j \in I_0$,
\eqn 
\yim\yjmp = t^{a_{ij}} \yjmp \yim + \delta_{ij}(1-t^{2}),
\eneqn
\item 
for $p > m+1$ and $i,j \in I_0$,
\eqn 
\yim \yjp = t^{(-1)^{p-m+1}a_{ij}} \yjp \yim.
\eneqn
\end{enumerate}

\begin{remark}
We change $t$ into $t^{-1}$ in the presentation in \cite{HL13}.
\end{remark}

\begin{theorem}[{\cite[Theorem 7.3]{HL13}}] \label{thm: braid group action}
Let $Q$ be a Q-data.
Then there is an isomorphism
$\iota_Q\col \KT(\g_0)\isoto \Kt(\g)$
such that
$\iota_Q(\yim)$ is equal to
$\qt{\D^m \F_Q(L(i))}$, where $L(i)$ is the simple module in $R_{\g_0}\smod$
corresponding to $i\in I_0$.
\end{theorem}

Let $B(\g_0)$ be the Braid group associated with $\g_0$. 
It is generated by $\set{\brd_i}{i\in I_0}$ with the defining relations
\eqn
\brd_i\brd_j\brd_i&=&\brd_j\brd_i\brd_j \quad \text{if $a_{ij}=-1$,} \\
\brd_i\brd_j&=&\brd_j\brd_i \quad \text{if $a_{ij}=0$.}
\eneqn
One of our main theorems is the following.
\begin{theorem} \label{thm: braid group action}
The Braid group $B(\g_0)$ acts on $\KT(\g_0)$ by the following formulas: 
\eqn 
\begin{aligned}
\brd_i(\yjm)=
&\begin{cases}
y_{j,m+\delta_{ij}}  & \text{if} \  a_{ij}\neq -1, \\[1ex]
\dfrac{t^{1/2} \yjm\yim-t^{-1/2}\yim\yjm}{t-t^{-1}}  & \text{if} \  a_{ij}=-1,
\end{cases}
\end{aligned}
\eneqn
\eqn 
\begin{aligned}
\brd_i^{-1}(\yjm)= 
&\begin{cases}
y_{j,m-\delta_{ij}}  & \text{if} \  a_{ij}\neq -1, \\[1ex]
\dfrac{t^{1/2} \yim\yjm-t^{-1/2}\yjm\yim}{t-t^{-1}}  & \text{if} \  a_{ij}=-1. 
\end{cases}
\end{aligned}
\eneqn
\end{theorem}

\Th\label{Saito}
Let $i$ be a sink of a Q-data $Q$.
Then the following diagrams commute:
\eqn
&&\xymatrix@C=5ex@R=4ex{
K(R_{\g_0}\gmod)\ar[rr]^-{j_Q}\ar[drr]_-{j_{s_iQ}}
&&\Kt(\g)\ar@{<-}[r]^\sim_{\iota_Q}&\KT(\g_0)
\ar[d]^{\brd_i}\\
&&\Kt(\g)\ar@{<-}[r]^\sim_{\iota_Q}&\KT(\g_0),
}\\
&&\xymatrix@C=5ex@R=1ex{
K({}_iR_{\g_0}\gmod)\;\ar[dd]^\bwr_{\mathbb{T}_i}\ar@{^{(}->}[r]&K(R_{\g_0}\gmod)\ar[dr]^-{j_Q}\\
&&\Kt(\g)\;.\\
K({}^iR_{\g_0}\gmod)\;\ar@{^{(}->}[r]&K(R_{\g_0}\gmod)\ar[ur]_-{j_{s_iQ}}
}
\eneqn
\enth
Here, ${}_iR_{\g_0}\gmod$ \ro resp.\  ${}^iR_{\g_0}\gmod$\rf\ is the full subcategory of
$R_{\g_0}\gmod$ consisting of graded modules $M$ with $E_i^*M=0$
\ro resp.\ $E_iM=0$\rf, and $\mathbb{T}_i$ is the reflection functor
due to S.~Kato \cite{Kato14,Kato} \ro cf.\ Y.~Saito \cite{Saito}\rf.
For $E_i$ and $E^*_i$, see for example, \cite{KKOP18}.

\section{The category $\T_N$ and reflection functors}
Let $J$ be the index set of simple roots of the root system $A_\infty$. One can identify $J$ with $\Z$ 
and the root lattice $\rootl$ is the subspace of  $\soplus_{a\in \Z} \Z\eps_a$ generated by $\al_a=\eps_a-\eps_{a+1}$ for $a\in\Z$.
 Let $R^{A_\infty}$ be the symmetric quiver Hecke algebra of type $A_\infty$ 
over $\cor$ with  the choice of parameters
\eqn 
Q_{ij}(u,v)=\delta(i\neq j)(u-v)^{\delta(j=i+1)}(v-u)^{\delta(j=i-1)}
\eneqn
for $i,j\in J$. It is a family $\{R^{A_\infty}(\beta)\}_{\beta \in \rootl^+}$ of associative $\Z$-graded $\cor$-algebras, where $\rootl^+=\sum_{i\in J} \Z_{\ge 0} \al_i$ is the positive root lattice of type $A_\infty$.
Each $R^{A_\infty}(\beta)$ is generated by $\{e(\nu)\}_{\nu \in J^\beta}$, $\{x_k\}_{1\le k \le n}$ and $\{\tau_m\}_{1\le m\le n-1}$, where 
$n=|\beta|\seteq\sum_{i\in I} n_i $ with $\beta=\sum_{i\in J} n_i \al_i$, and    $J^\beta\seteq\set{\nu \in  J^n}{\al_{\nu_1}+\cdots +\al_{\nu_n}=\beta}$. 
See \cite{KKK18A} for a set of defining relations of  $R^{A_\infty}(\beta)$.
Note that there is an embedding of $R^{A_\infty}(\beta)\tens R^{A_\infty}(\gamma)$  into $R^{A_\infty}(\beta+\gamma)$. 
Hence the category $\catA=\bigoplus_{\beta \in \rootl^+} R^{A_\infty}(\beta)\gmod $ of finite-dimensional graded $R^{A_\infty}$-modules is a monoidal category whose tensor product is given by the \emph{convolution product}: 
$$M \conv N\seteq R^{A_\infty}(\beta+\gamma) \tens_{R^{A_\infty}(\beta) \tens R^{A_\infty}(\gamma)} \left( M\tens N \right).$$
For $M \in R^{A_\infty}(\beta)\gmod$, we set $\wt(M)\seteq-\beta$. 

For each pair of integers $a,b $  with $a\le b$, let $[a,b]$ be the interval $\set{k\in \Z}{ a \le k\le b}$, and call it a \emph{segment}.
For each segment $[a,b]$, let $L(a,b)$ be
  the one-dimensional simple graded $R^{A_\infty}$-module generated by a vector $u(a,b)$  such that $e(\nu) u(a,b) = \delta(\nu=(a,\ldots,b))u(a,b)$. 
We set $L(a)\seteq L(a,a)$ for $a \in \Z$. 
For each $N\ge 2$, let $\mathcal S_N$ be the smallest subcategory of $\catA$ which is 
stable under taking convolution, subquotients, extensions, and containing $\set{L(a,b)}{b-a+1 > N}$.
Then the quotient category $\catA  / \mathcal S_N$ equips with a new tensor product $\star$ 
 given by
 \eqn
 X \star Y \seteq t^{B(\wt(X),\wt(Y))}X\conv Y,
 \eneqn
 where
   $B(x,y) \seteq -\sum_{k >0} (S^k x,y)$ for $x,y\in \soplus_{a\in \Z} \Z\eps_a$ and 
   $S$ is an automorphism on $\soplus_{a\in \Z} \Z\eps_a$ given by $S(\eps_a) \seteq\eps_{a+N}$.
The category $\T_N$ is constructed   in \cite{KKK18A} as   a  localization of the monoidal category $(\catA  / \mathcal S_N ,\star)$. 
The objects of $\T_N$ is the same with the ones of $\catA  / \mathcal S_N$. 
The group of  morphisms is  given by 
$$\Hom_{\mathcal T_N}(X,Y) \seteq\indlim_{\la,\mu} \Hom_{\catA / \mathcal S_N}(X\conv P^\la, Y\conv P^\mu),$$
where 
$P^\nu\seteq\conv_{a\in \Z} L(a,a+N-1)^{\circ \nu_a} $ for $\nu \in (\Z_{\ge 0})^{\oplus J} $ and  the limit runs over  all the pairs $(\la,\mu)$ such that $\wt(X \conv P^\la) =\wt(Y \conv P^\mu)$.
It turns out that  $\T_N$ is an abelian rigid monoidal category with a tensor product $\star$. 
We denote the right dual (resp. left dual) of $X$ by $\D(X)$ (resp. $\D^{-1}(X)$).
Note that  $L(a,a+N-1)\simeq \one$ in $\T_N$ for all $a\in \Z$.
 We have a chain of functors
$$\catA \To[\mathcal Q_N] \catA/\mathcal S_N \To[\Upsilon_N]  \T_N.$$ 
The composition will be denoted by $\Omega_N$.
Note that  the Grothendieck ring $K(\T_N)$  is a $\Z[t^{\pm1}]$-algebra on which $t$ acts by the grading shift.

From now on, let $\g$ be the affine Kac-Moody algebra of type ${A_{N-1}^{(1)}}$.
We regard $\T_N$ as a $\Z$-graded lifting of $\catC_\g^{0}$ as a rigid monoidal category. 
Indeed there exists a monoidal functor $\F_N\col\T_N\to \catC_\g^0$ which sends simples to simples.
It induces an isomorphism of $\C(t^{1/2})$-algebras 
$[\F_N] \col\C(t^{1/2})\tens_{\Z[t^{\pm1}]}K(\T_N)\isoto\Kt(\g)$ (\cite[Theorem 4.33]{KKK18A}).  Under the isomorphism, the generator $y_{i,m}$ corresponds to $[\D^mL(i)]$ for $i\in \Z, \ m\in \Z$.

 For a pair $(M,N)$ of objects in an abelian monoidal category we denote by $M\hconv N$ the head of $M\tens N$ and by $M\sconv N$ the socle of $M\tens N$, respectively.

We show that there is a family of autofunctors on $\T_N$ which recover the braid group action on the quantum Grothendieck ring $\Kt(\g)$. 
 For this purpose, we adjoin a  formal object $t^{1/2} \one$ into $\T_N$ such that $t^{1/2} \one \star t^{1/2} \one \simeq t \one$. Then the \emph{grading shift by $1/2$ of $X$} is given by $X \to t^{1/2}\one \star X$.

\begin{theorem} \label{thm: reflection functors}
For $i \in \Z$, there exists a monoidal functor 
$\Refl_i \col \T_N \to \T_N$ satisfying
\eqn
&&\Refl_i (L(j)) 
\simeq \begin{cases}
\D L(j)& \text{if} \ j \equiv i\ \mathrm{mod}\; N, \\
 t^{1/2} \left( L(j\mp 1) \hconv L(j)\right)& \text{if} \ j \equiv i\pm 1\ 
\mathrm{mod}\; N, \\
L(j)& \text{otherwise.} \\
\end{cases}
\eneqn

The functor $\Refl_i$ has an inverse
$\Reflinv_i \col \T_N \to \T_N$ satisfying
\eqn
&&\Reflinv_i (L(j))
\simeq\begin{cases}
\D^{-1} L(j)& \text{if} \ j \equiv i \ \mathrm{mod}\; N, \\
 t^{1/2}  \left( L(j) \hconv L(j\mp1) \right)& \text{if} \ j \equiv j \pm 1
\ \mathrm{mod}\;  N, \\
L(j)& \text{otherwise.} \\
\end{cases}
\eneqn
\end{theorem}

Let us explain briefly how to construct the functors $\Refl_i$.
For each $j\in J$, denote $\bar M_j$  the $R^{A_\infty}$-module
$t^{-1} L(j+1,j+N-1)$ if $j\equiv i \mod N$, 
 $t^{1/2} \left( L(j\mp 1) \hconv L(j)\right)$ if $j \equiv i\pm 1 \mod N$ and 
$L(j)$ otherwise.
For each $\beta \in \rootl^+$ and $\mu=(\mu_1,\ldots,\mu_m) \in J^\beta$, set
\eqn
\Delta(\mu) =M_{\mu_1}\conv \cdots \conv M_{\mu_m}, \ \text{and} \
\Delta(\beta) =\soplus_{\mu \in J^\beta} \Delta(\mu),
\eneqn
where 
$M_j$ is the affinization of $\bar M_j$.
Then along a similar line with \cite[Section 4]{KP},  one can show that there exists a ring homomorphism
\eqn
(R^{\A_\infty}(\beta))^\opp \to \End_{\catA^{\mathrm big} /\mathcal S^{\mathrm big}_N} (\mathcal Q_N(\Delta(\beta))),
\eneqn
where $\catA^{\mathrm big} /\mathcal S^{\mathrm big}_N$ is an infinite analogue of $\catA/\mathcal S_N$.
Let $\mathcal R'_\beta \col  R^{A_{\infty}}(\beta) \gmod \to  \catA /\mathcal S_N$ be the restriction of a left adjoint of the functor $\Hom_{\catA{\mathrm big}/\mathcal S{\mathrm big}_N} (\mathcal Q_N(\Delta(\beta)), -) $.
Then we obtain a monoidal functor  $\mathcal R\col \catA \to \T_N $, the composition $$\catA \To[\bigoplus_{\beta \in \rootl^+}\mathcal R'_\beta] \catA /\mathcal S_N \To[\Upsilon_N] \T_N.$$
Note that   the family 
$\{\bar M_j\}_{j\in J}$ of objects in $\T_N$ satisfies  for any $a\in J$  that (1) $\bar M_a \star \bar M_{a+1} \star  \cdots \star \bar M_{a+N-1} \simeq \one$, (2)  $\hd(\bar M_a \star \bar M_{a+1} \star  \cdots \star \bar M_{a+k-1}) \star \bar M_{a+k}$ is not simple for $1\le k\le N-1$, and (3) $\D^2(\bar M_{a}) \simeq \bar M_{a+N}$. 
A similar argument as the one in \cite[Section 6.1]{KKOP19B} shows that there is a monoidal functor $\Refl_i\col\T_N \to \T_N$ such that  $\mathcal R \simeq \Refl_i \circ \Omega_N$. 
\medskip

Recall that there is an automorphism  $\shift\col\T_N \to \T_N$ given by $L(j) \mapsto L(j+1)$ for all $j\in \Z$. It satisfies that  $\shift^N \simeq \D^2$.
The functors $\{ \Refl_i   \ | \ i\in \Z\}$ satisfy the following properties.
\begin{prop} We have
\bnum
\item $\Refl_{i+1} \simeq \shift \circ \Refl_i \circ \shift^{-1}$ for $i\in \Z$,
\item $\Refl_{i} \circ \D  \simeq\D \circ \Refl_i$ for $i\in \Z$,
\item $\Refl_i \simeq \Refl_{N+i}$ for $i\in \Z$,
\item $\Refl_{1}\Refl_2\cdots \Refl_{N-1}\simeq\shift$,
\item $\Refl_i \circ \Refl_j \simeq \Refl_j \circ \Refl_i $ for $|i-j|>2$,
\item $\Refl_i \circ \Refl_{i+1} \circ \Refl_i \simeq \Refl_{i+1} \circ \Refl_i \circ \Refl_{i+1}$ for $i\in \Z$.
\end{enumerate}
\end{prop}

\medskip
The family of functors $\{\Refl_i\}_{1\le \i \le N-1}$ recovers the braid group action in Theorem \ref{thm: braid group action} in the case of type $A_{N-1}$.
\begin{theorem} \label{thm:functors recover the action}
For each $1\le i \le N-1$  the $\Z[t^{\pm 1/2}]$-algebra automorphism on $K(\T_N)$ induced by  $\Refl_i$ is equal to $\brd_i$ in 
{\rm Theorem \ref{thm: braid group action}} under the isomorphism 
$$[\F_N] \col  \C(t^{1/2})\tens_{\Z[t^{\pm1}]}K(\T_N)\isoto\Kt(\g).$$
\end{theorem}

\medskip
\section{Reflections by root modules}
Recall that for each pair of non-zero modules $(X,Y)$ of $\catA$, there exists a distinguished nonzero morphism $\rmat{X,Y}\col t^{\Lambda(X,Y)}X \conv Y \to Y \conv X$ called the \emph{r-matrix} (\cite{KKK18A}). 
Here, $t$ is the grading shift functor.
We have 
${\Omega_N(\rmat{X,Y}}) \col t^{\Lambda_N(X,Y)}X \star Y \to Y \star X$ in $\T_N$,
where  $\Lambda_N(X,Y) = \Lambda(X,Y)-B(\wt(X),\wt(Y)) + B(\wt(Y),\wt(X))$.

For a pair $(X,Y)$ of objects in $\T_N$, set $$\de(X,Y)\seteq\dfrac{1}{2}(\Lambda_N(X,Y) + \Lambda_N(Y,X)).$$
Note that $\de(X,Y)= \dfrac{1}{2}(\Lambda(X,Y) + \Lambda(Y,X))$ if $\Omega(\rmat{X,Y}) \neq 0$.

A simple object $X$ in an abelian monoidal category is called \emph{real} if $X\tens X$ is  simple.
A real simple object $L$ in $\T_N$ is called a \emph{root module} if 
\eqn 
\de(L, \D^k(L)) =\delta(k =\pm1).
\eneqn
For example, the objects $L(a,b)$ with $b-a+1 < N$ are  root modules. 
If $L$ is a root module, then $\D(L)$, $\D^{-1}(L)$ and $\Refl_i(L)$ for $i\in \Z$ are  root modules.

The following is the main theorem of this section.
\begin{theorem} \label{thm:S_i(X)}
Let $X$ be a simple object in $\T_N$.
For $i\in \Z$, if 
\eqn \de(\D^k(L(i)), X) =n \delta(k=a)\eneqn
for some $n\ge 0$ and $a \in \Z$, then
$$\Refl_i (X) \simeq (\D^a L(i))^{\circ n} \hconv X$$
up to a multiple of a power of $t$.
\end{theorem}

The following is 
one of  the applications of Theorem \ref{thm:S_i(X)}.
\begin{theorem} \label{thm: reflection by root modules}
Let $[L]$ belongs to the orbit of $L(i)$ for some $1\le i\le N-1$
under the braid group $B(A_{N-1})$ action in Theorem~\ref{thm: braid group action}. 
Then there is an automorphism $s_L$ on $K(\T_N)$ such that
\bnum
\item $s_{L(i)} = s_i$ for $1\le i\le N-1$.
\item $s_{(s_{L}(L'))} =s_L \circ s_{L'} \circ s_{L}^{-1}$ if $L'$ also satisfies the condition in the theorem.
\item $s_{\D^a L} =s_{L}$ for all $a \in \Z$
\item $s_{L}([X]) = [(\D^{a} L )^{\circ n} \hconv X]$ up to a power of $t$, 
if $\de(\D^k L, X) = n\delta(k=a)$ for some $n\ge 0$ and $a \in \Z$.
\end{enumerate}
\end{theorem}

\section{Conjectures}

Let $\uqpg$ be an arbitrary quantum affine algebra.
We say that a real simple module $L$ in $\catC^0_\g$ is a {\em root module} if
$\de(\dual^kM,M)=\delta(k=\pm1)$ for any $k$.

\Conj
For any root module $L\in\catC^0_\g$, there exists a monoidal autofunctor
$\Refl_L$ of $\cat$ which satisfies the following conditions:
\bna
\item $\Refl_L$ satisfies similar properties in
Theorem~\ref{thm:S_i(X)} and Theorem~\ref{thm: reflection by root modules}.
\item {\rm(Braid relation)}
For root modules $L$ and $L'$,
\be[{\rm(1)}]
\item 
if $\de(\dual^kL,L')=0$ for any $k\in \Z$, then
$$\Refl_L\circ\Refl_{L'}\simeq\Refl_{L'}\circ\Refl_{L},$$
\item
if $\de(\dual^kL,L')=\delta(k=0)$ for any $k\in \Z$, then
$$\Refl_L\circ\Refl_{L'}\circ\Refl_L\simeq\Refl_{L'}\circ\Refl_{L}\circ\Refl_{L'}.$$
\ee
\item 
Let $Q$ be a Q-data, and $L\seteq\F_Q(L(i))$.
Then the automorphism of $\Kt(\g)$ induced by $\Refl_L$ coincides with
$\brd_i$, i.e.,
the following diagram commutes:
$$\xymatrix{
\KT(\g_0)\ar[r]^\sim_{\iota_Q}\ar[d]_{\brd_i}&\Kt(\g)\ar[d]_{\Refl_L}\\
\KT(\g_0)\ar[r]^\sim_{\iota_Q}&\Kt(\g).
}$$
\ee
\enconj

\end{document}